# EXISTENCE AND CONSTRUCTION OF RANDOMIZATION DEFINING CONTRAST SUBSPACES FOR REGULAR FACTORIAL DESIGNS

By Pritam Ranjan, Derek R. Bingham and Angela M. Dean

*Acadia University, Simon Fraser University and Ohio State University*

Regular factorial designs with randomization restrictions are widely used in practice. This paper provides a unified approach to the construction of such designs using *randomization defining contrast subspaces* for the representation of randomization restrictions. We use finite projective geometry to determine the existence of designs with the required structure and develop a systematic approach for their construction. An attractive feature is that commonly used factorial designs with randomization restrictions are special cases of this general representation. Issues related to the use of these designs for particular factorial experiments are also addressed.

**1. Introduction.** The use of factorial experiments in situations which require randomization restrictions (e.g., block designs or split-plot designs) has been illustrated in the literature for many years [see, e.g., Cochran and Cox (1957), Chapters 6–8 and Addelman (1964)]. More recently, the construction of designs under various optimality and efficiency criteria has been discussed for blocked factorial and fractional factorial designs by Bisgaard (1994), Sitter, Chen and Feder (1997), Sun, Wu and Chen (1997), Chen and Cheng (1999) and Cheng, Li and Ye (2004); for factorial and fractional factorial split-plot designs by Huang, Chen and Voelkel (1998), Bingham and Sitter (1999) and Stapleton, Lewis and Dean (2009); for strip-plot designs by Miller (1997); and for split-lot designs by Mee and Bates (1998) and Butler (2004). Although these various designs maintain the same factorial treatment structure, their randomization structures are different.

Bingham et al. (2008) proposed the use of *randomization defining contrast subgroups* (RDCSGs) to describe the randomization structure of factorial designs. These can be viewed as generalizations of block defining contrast subgroups [see, e.g., Sun, Wu and Chen (1997)]. In this article, we propose a projective geometric generalization of the RDCSGs, which we refer to as *randomization defining contrast subspaces* (RDCSSs).

For complicated randomization structures, the existence of desirable designs can be difficult to determine; for example, Bingham et al. (2008) were forced to search over the set of all possible designs to find those satisfying the required restrictions. Here, we develop theoretical results for the existence of factorial designs with given randomization restrictions within a unified framework, and we provide a direct method of construction of such designs. The theory for the existence and construction of such designs is given in terms of full factorial designs, but we show in Section 6 that the results apply equally well to fractional factorial designs.

In the next section, a brief description of RDCSSs is presented. Section 3 considers the impact of the randomization restrictions on data analysis and motivates the preference for nonoverlapping RDCSSs over the various stages of randomization. In Section 4, theoretical results for the existence of designs with randomization restrictions are developed. The theoretical framework is constructive, thereby allowing for the identification of designs in practical settings. Methods of construction of such designs for specific factorial experiments are developed in Section 5. The article finishes with a discussion of fractional factorial designs in Section 6 and concluding remarks in Section 7.

**2. Randomization-defining contrast subspaces.** A natural design choice for a *p*-factor experiment is a completely randomized design with trials performed in a random order. However, this is not always feasible as, for example, when the experimental units are not homogeneous (as in block designs), the levels of some factors cannot be changed as rapidly as that of others (as in split-plot designs) or when subsets of factors must be held fixed at different stages of experimentation (as in split-lot designs).

Several different approaches have been taken in the construction of factorial designs with randomization restrictions. For example, blocking in factorial and fractional factorial experiments have been studied using Abelian groups and vector spaces [see, e.g., Dean and John (1975), Bailey (1977) and Voss and Dean (1987)], and using finite geometries [see, e.g., Bose (1947), Srivastava (1987) and Mukerjee and Wu (1999)]; the construction of strip-plot designs via latin square fractions has been discussed by Miller (1997); split-plot designs for factorial treatments have been constructed using computer search and group theory by Addelman (1964), Bingham and Sitter (1999), Bisgaard (2000), Huang, Chen and Voelkel (1998) and Stapleton, Lewis and Dean (2009). For the construction of split-lot designs, linear graphs were used by Taguchi (1987), combinations of linear graphs and cyclic groups by Mee and Bates ([1998](#)) and a grid-representation technique by Butler ([2004](#)). In this paper, we use finite projective geometry to study $q^p$ factorial designs with randomization restrictions for prime or prime power $q$.

As in Dey and Mukerjee ([1999](#)), Chapter 8, let $b$ be a $p$-dimensional pencil over the Galois field $GF(q)$. For $\alpha(\neq 0) \in GF(q)$, $b$ and $\alpha b$ represent the same pencil carrying $q - 1$ degrees of freedom. A pencil $b$ represents an $r$-factor interaction if $b$ has exactly $r$ nonzero elements. The set of all $p$-dimensional pencils over $GF(q)$ forms a $(p - 1)$-dimensional finite projective geometry and is denoted by $PG(p - 1, q)$. There are $(q^p - 1)/(q - 1)$ points in $PG(p - 1, q)$, where points correspond to pencils.

Since two-level factorial designs are the most common in practice, we focus on $q = 2$, although most of the results in this article hold for a prime or prime power $q$. For $q = 2$, a pencil $b$ with $r$ nonzero elements corresponds to a unique $r$-factor interaction in a $2^p$ factorial design with a single degree of freedom. Thus, the set of all effects (excluding the grand mean) of a $2^p$ factorial design is equivalent to the set of points in $PG(p - 1, 2)$ and will be referred to as the *effect space* $\mathcal{P}$.

The restrictions on the randomization of experimental runs is equivalent to the grouping of experimental units into sets of trials. We consider the usual approach of using independent effects from $\mathcal{P}$ to define the groupings. Blocked factorial designs, for example, use $2^t$ ($t < p$) combinations of $t$ independent effects from $\mathcal{P}$ to divide $2^p$ treatment combinations into $2^t$ blocks. These factorial effects are then completely confounded with block effects and represent $t$ randomization restriction "factors." The set $S$ of all nonnull linear combinations of these $t$ randomization restriction factors in $\mathcal{P}$, over $GF(2)$, forms a $(t - 1)$-dimensional projective subspace of $\mathcal{P} = PG(p - 1, 2)$, which we call a *randomization defining contrast subspace* (RDCSS).

EXAMPLE 1.   Consider an experiment arranged as a $2^5$ factorial split-plot design, where $A$, $B$ are whole-plot factors and $C$, $D$, $E$ are sub-plot factors. The effect space is $\mathcal{P} = \langle A, B, C, D, E \rangle$ and the RDCSS that imposes the randomization restrictions is $S = \langle A, B \rangle = \{A, B, AB\}$, where the notation $\langle a_1, \ldots, a_k \rangle$ denotes the projective space spanned by $a_1, \ldots, a_k$ (or, for two-level factorial designs, the set of all interactions of $a_1, \ldots, a_k$). The notation $A$, $AB$, and so on, represent the main effect of $A$, the interaction of $A$ and $B$, and so on. Since there are $t = 2$ whole-plot factors, the set of all experimental units is partitioned into $2^t = 4$ subsets (batches, whole-plots, etc.), and each subset consists of $2^5/2^2 = 8$ experimental units. The four subsets, say $\mathcal{B}_1$, $\mathcal{B}_2$, $\mathcal{B}_3$ and $\mathcal{B}_4$, consist of experimental units corresponding to $(\theta_A(i), \theta_B(i)) = (0,0), (0,1), (1,0)$ and $(1,1)$, respectively, where, $\theta_\delta(i)$ is the $i$th row of the column corresponding to the factorial effect $\delta$ in the model matrix $X$ [see, e.g., Bingham et al. ([2008](#)) and page 26 in Ranjan ([2007](#))].

There may be more than one stage of randomization restriction in a factorial experiment where the randomization structure can be characterized by its RDCSSs. For a $2^p$ factorial design with $m$ stages of randomization, the $m$ RDCSSs can be denoted by the projective subspaces $S_1, \ldots, S_m$ contained in the effect space $\mathcal{P}$. For each $i = 1, \ldots, m$, the size of $S_i$ is $2^{t_i} - 1$ with $0 < t_i < p$. Then, at stage $i$, the experimental units are partitioned into $|S_i| + 1$ sets (e.g., batches, whole-plots and blocks) due to $S_i$, where the size of each set is $(|\mathcal{P}| + 1)/(|S_i| + 1) = 2^{p-t_i}$.

EXAMPLE 2. Consider a $2^5$ factorial experiment with randomization structure defined by a strip-plot design [Miller ([1997](#))], where the row configurations are defined by a $2^2$ design in factors $A$, $B$ and the column configurations by a $2^3$ design in factors $C$, $D$, $E$. In this setting, the two RDCSSs are $S_1 = \langle A, B \rangle$ and $S_2 = \langle C, D, E \rangle$, and the effect space is $\mathcal{P} = \langle A, B, C, D, E \rangle$.

Although the treatment structure for Examples 1 and 2 are the same, the randomization restrictions induce different error structures, thereby having an impact on the analysis [see, e.g., Milliken and Johnson ([1984](#)), Chapter 4].

**3. Modeling issues.** In this section, we show how the distribution of the least squares estimator of a factorial effect is related to the RDCSSs. This motivates a design strategy and, in particular, suggests a preference for nonoverlapping RDCSSs.

3.1. *Model.* Consider a single-replicate $2^p$ factorial design with the linear regression model as the response model of interest; that is,

$$(1) \qquad Y = X\beta + \varepsilon,$$

where $X$ denotes the $n \times 2^p$ model matrix and $\beta = (\beta_0, \beta_1, \ldots, \beta_{2^p-1})'$ is the $2^p \times 1$ vector of parameters corresponding to the grand mean, the factorial main effects and interactions. In this section, the factor levels 0 and 1 are recoded as $+1$ and $-1$, respectively. For a single-replicate $2^p$ factorial experiment, $n = 2^p$, and the model matrix, $X$, is a Hadamard matrix which satisfies $X'X = nI_n$, where $I_n$ is an $n \times n$ identity matrix. Without loss of generality, any $p$ independent columns of $X$ can be selected to represent the main effect contrasts of the $p$ factors [see, e.g., Dean and Voss ([1999](#)), Section 15.6]. On rearranging the columns of $X$, let $X$ be $\{c_0, c_1, \ldots, c_p, c_{p+1}, \ldots, c_{n-1}\}$, where the column vector, $c_0$, consists of all 1's corresponding to the grand mean, columns $c_1, \ldots, c_p$ refer to the independent main effect contrasts of the $p$ factors and the remaining columns, $c_{p+1}, \ldots, c_{n-1}$, represent the interaction contrasts obtained as element-wise products of subsets of $c_1, \ldots, c_p$. The vector $Y$ denotes the vector of response variables and $\varepsilon$ the error vector, whose distribution is discussed below.

For a factorial design with $m$ stages of randomization and RDCSSs denoted by $S_i$, $i = 1, \ldots, m$, the error vector $\varepsilon$ can be written as a sum of $m + 1$ independent error vectors, $\varepsilon = \varepsilon_0 + \varepsilon_1 + \cdots + \varepsilon_m$. The $n \times 1$ vector $\varepsilon_0$ denotes the replication error vector, and $\varepsilon_i$ $(1 \leq i \leq m)$ is the error vector associated with the randomization restriction characterized by $S_i$. Since $S_i$ creates a partition of the $n$ experimental units into $2^{t_i}$ batches (blocks, whole-plots, etc.), the error vector $\varepsilon_i$ can be written as $N_i \epsilon_i$, where $\epsilon_i$ is a $2^{t_i} \times 1$ vector of errors associated with each of the $2^{t_i}$ batches; the $n \times 2^{t_i}$ matrix $N_i$ is called the $i$th *incidence matrix* and has elements defined as follows:

$$(2) \qquad \begin{aligned} (N_i)_{rl} &= 1, && \text{if the } r\text{th experimental unit belongs to the } l\text{th batch} \\ && & \text{at the } i\text{th stage of randomization,} \\ &= 0, && \text{otherwise,} \end{aligned}$$

for $r = 1, \ldots, n$ and $l = 1, \ldots, 2^{t_i}$. Thus, $\varepsilon$ can be written as

$$(3) \qquad \varepsilon = \varepsilon_0 + \varepsilon_1 + \cdots + \varepsilon_m = \epsilon_0 + N_1 \epsilon_1 + \cdots + N_m \epsilon_m,$$

where $\epsilon_0, \epsilon_1, \ldots, \epsilon_m$ are independent, and $\epsilon_0 = \varepsilon_0 \sim N(0_n, I_n \sigma^2)$, where $0_n$ is an $n \times 1$ vector of zeros and $I_n$ is an $n \times n$ identity matrix. For $i = 1, \ldots, m$, we assume that $\epsilon_i \sim N(0_{2^{t_i}}, I_{2^{t_i}} \sigma_i^2)$, where $0_{2^{t_i}}$ is a vector of $2^{t_i}$ zero elements. It follows that $\varepsilon \sim N(0_n, \Sigma_y)$, where

$$(4) \qquad \Sigma_y = \sigma^2 I_n + \sum_{i=1}^{m} \sigma_i^2 N_i N_i'.$$

Lemma 1 gives some properties of the incidence matrices which are needed for subsequent results. The proof of Lemma 1 is straightforward and omitted.

LEMMA 1.   *Consider a $2^p$ factorial design with $n = 2^p$ runs and $m$ levels of randomization restrictions defined by RDCSSs $S_1, \ldots, S_m$. Let $N_1, \ldots, N_m$ be the incidence matrices corresponding to $S_i$'s as defined in (2). Then, $N_i' N_i = n_i I_{2^{t_i}}$, where $n_i = 2^{p - t_i}$ is the number of 1's in each column of $N_i$, and $I_{2^{t_i}}$ is a $2^{t_i} \times 2^{t_i}$ identity matrix.*

3.2. *Distribution of effect estimators.* The most natural way to estimate the factorial effect parameters is to use the *generalized least squares* (GLS) estimator $\hat{\beta} = (X' \Sigma_y^{-1} X)^{-1} X' \Sigma_y^{-1} Y$. Due to the assumptions on the error vectors, $\hat{\beta} \sim N(\beta, \mathrm{Var}(\hat{\beta}))$. For a single-replicate $2^p$ factorial design, the GLS estimator coincides with the ordinary least squares (OLS) estimator $\tilde{\beta} =$ $(X'X)^{-1}X'Y$. This follows since $X$ is an $n \times n$ matrix and, consequently, $X'X = nI$ implies that $X' = nX^{-1}$. Then, $\tilde{\beta} = n^{-1}X'Y$, and

$$\hat{\beta} = (X'\Sigma_y^{-1}X)^{-1}X'\Sigma_y^{-1}Y = X^{-1}\Sigma_y(X')^{-1}X'\Sigma_y^{-1}Y = n^{-1}X'Y = \tilde{\beta}.$$

Thus, the variance–covariance matrix of the factorial effect estimator $\hat{\beta}$ is $\mathrm{Var}(\hat{\beta}) = \mathrm{Var}(\tilde{\beta}) = (X'\Sigma_y X)/n^2$, and

$$\hat{\beta} \sim N(\beta, X'\Sigma_y X/n^2), \tag{5}$$

where $\Sigma_y$ is defined in (4).

If $S_i$'s are projective subspaces in $\mathcal{P}$, it is sometimes possible to select the RDCSSs $S_1, \ldots, S_m$ so that $S_{ij} = S_i \cap S_j = \phi$ for all $i \neq j$. It will become clear later that these cases are of specific interest to practitioners since, in this case, all the effects in $S_i$ have variances that are functions of $\sigma_i^2$ only. However, when this condition does not hold (see Theorem 2, below), the variance of the effects in $S_{ij}$ are impacted by both $\sigma_i^2$ and $\sigma_j^2$.

THEOREM 1. *In a single-replicate $2^p$ factorial design, let the randomization restrictions be defined by $S_1, \ldots, S_m$ contained in $\mathcal{P}$. Then, for any two effects $E_j$ and $E_k$ in the effect space $\mathcal{P}$, the corresponding parameter estimators $\hat{\beta}_{E_j}$ and $\hat{\beta}_{E_k}$ have independent normal distributions.*

PROOF. Since $\hat{\beta}$ has a multivariate normal distribution, it is enough to show that $\mathrm{cov}(\hat{\beta}_{E_j}, \hat{\beta}_{E_k}) = 0$. Let $E_j$ denote the factorial effect corresponding to $j$th column $c_j$ $(j \geq 1)$ of $X$, and let $S_1, \ldots, S_m$ be subspaces in $\mathcal{P}$. If column $c_j$ of $X$ is orthogonal to the randomization restrictions in $S_i$, then $c_j'N_i = 0$. Whereas, for each $i = 1, \ldots, m$, if $E_j \in S_i$, then, as in Dean (1978), $c_j$ may be expressed as $c_j = N_i a_j^{(i)}$, where $a_j^{(i)}$ is a $2^{t_i} \times 1$ vector, and $a_j^{(i)'}\mathbf{1} = 0$ with $\mathbf{1}$ a vector of 1's. Also, since $c_j'c_k = 0$ for all $j \neq k$, it follows from Lemma 1 that $a_j^{(i)'}a_k^{(i)} = 0$ for all $j \neq k$ and $i = 1, \ldots, m$.

Define the index set $T_{jk} = \{i : 1 \leq i \leq m, E_j \in S_i \text{ and } E_k \in S_i\}$. Then, using the equivalence of the GLS and OLS estimators for $\beta$ together with (4), (5) and Lemma 1, we have

$$n^2\,\mathrm{Cov}(\hat{\beta}_{E_j}, \hat{\beta}_{E_k}) = \sigma^2 c_j'c_k + \sum_{i=1}^{m}\sigma_i^2 c_j'N_iN_i'c_k$$

$$= \sum_{i \in T_{jk}}\sigma_i^2 a_j^{(i)'}N_i'N_iN_i'N_ia_k^{(i)} = \sum_{T_{jk}}\sigma_i^2 n_i^2 a_j^{(i)'}a_k^{(i)} = 0$$

and the proof follows. □

Next, the main result of this section establishes the link between the variances of the effect estimators for a $2^p$ full factorial design and the RDCSSs.

THEOREM 2.  *In a single-replicate $2^p$ factorial design, let the randomization restrictions be defined by $S_1, \ldots, S_m$ in $\mathcal{P}$. For an effect $E \in \mathcal{P}$, define an index set $\{T_E, E \in \mathcal{P}\}$ such that $T_E = \{i : 1 \le i \le m, E \in S_i\}$. Then,*

$$\operatorname{Var}(\hat{\beta}_E) = \begin{cases} \dfrac{\sigma^2}{n} + \displaystyle\sum_{i \in T_E} \dfrac{n_i}{n}\sigma_i^2, & \text{if } E \in \{S_1 \cup \cdots \cup S_m\}, \\[2ex] \dfrac{\sigma^2}{n}, & \text{if } E \in \mathcal{P} \setminus \{S_1 \cup \cdots S_m\}, \end{cases}$$

*where $\sigma^2$ and $\sigma_i^2$ are the variances corresponding to $\varepsilon_0$ and $\varepsilon_i$, respectively.*

PROOF.  Let $E$ be the effect corresponding to column $c$ in $X$. Following arguments similar to those of the proof of Theorem 1, if $E \in S_i$, we have

$$(6) \qquad c'c = a^{(i)'} N_i' N_i a^{(i)} = n_i a^{(i)'} a^{(i)} = n$$

and, if $E \notin S_i$, then $c' N_i = 0$. Thus,

$$n^2 \operatorname{Var}(\hat{\beta}_E) = n\sigma^2 + \sum_{i \in T_E} \sigma_i^2 a^{(i)'} N_i' N_i N_i' N_i a^{(i)}$$

$$= n\sigma^2 + \sum_{i \in T_E} \sigma_i^2 n \cdot n_i \qquad [\text{using Lemma 1 and (6)}]$$

as required. On the other hand, if $E \in \mathcal{P} \setminus \{S_1 \cup \cdots \cup S_m\}$, $T_E$ is the empty set and $\operatorname{Var}(\hat{\beta}_E) = \sigma^2/n$.  $\square$

We note that the results in Theorems 1 and 2 can be extended to regular designs with factors having $q$ levels for prime or prime power $q$. This is done by including in the model matrix a set of orthogonal contrasts for each effect, with each contrast scaled so that $c'c = n$.

A common strategy for the analysis of unreplicated factorial designs is the use of half-normal plots [Daniel (1959)]. These require that effects whose estimates appear on the same plot must be independent and have the same variance. By Theorems 1 and 2, when $S_{ij} = S_i \cap S_j = \phi, \forall i \neq j$, $m$ separate half-normal plots (plus an additional plot if $\mathcal{P} \setminus \{\bigcup_{i=1}^m S_i\} \neq \phi$) can be constructed to assess the significance of the effects. On the other hand, if $S_{ij} \neq \phi$ for some $i, j$, then the effects in $S_{ij}$ will have variances that are linear combinations of $\sigma_i^2$ and $\sigma_j^2$, creating a larger number of smaller sets of effects having the same variance. In this case, assessment of the significance of effects in $S_{ij}$ may have to be sacrificed due to lack of degrees of freedom [Schoen (1999)]. Thus, designs with disjoint RDCSSs ($S_{ij} = \phi$ for all $i \neq j$) are preferred. From a practical standpoint, determining whether or not a design with nonoverlapping RDCSSs exists is challenging [see, e.g., Bingham et al. (2008)] and actually finding such a design when it exists can be difficult. In the next section, we develop conditions for the existence of such designs.

**4. Existence of RDCSS.** In Section 4.1, results are given that focus on the maximum number of disjoint subspaces of equal size that are contained in the effect space $\mathcal{P} = PG(p-1, 2)$. New results are developed in Section 4.2 for the more general setting of unequal sized RDCSSs. This latter case is important in multistage experiments where the number of units $(n_i)$ in a batch at stage $i$ is not the same for all $i = 1, \ldots, m$. The results are constructive and thus allow experimenters to find designs in practice.

4.1. *RDCSSs, spreads and disjoint subspaces.* In most applications, the number of stages, $m$, of randomization is pre-specified by the experimenter. Thus, if a set, $\mathcal{S}$, of disjoint subspaces can be obtained, with $|\mathcal{S}| > m$, an appropriate subset of $\mathcal{S}$ can be chosen to construct the RDCSSs.

DEFINITION 1. For $1 \leq t \leq p$, a $(t-1)$-spread of the effect space $\mathcal{P}$ is a set, $\mathcal{S}$, of $(t-1)$-dimensional subspaces of $\mathcal{P}$, which partitions $\mathcal{P}$.

In a $(t-1)$-spread, $\mathcal{S}$, every element of $\mathcal{P}$ is contained in exactly one of the $(t-1)$-dimensional subspaces of $\mathcal{S}$. A $(t-1)$-spread is said to be nontrivial if $t > 1$. Given that a $(t-1)$-spread, $\mathcal{S}$, of $\mathcal{P}$ exists, the size of $\mathcal{S}$ is $|\mathcal{S}| = (2^p - 1)/(2^t - 1)$. A necessary and sufficient condition for the existence of a $(t-1)$-spread is that $t$ divides $p$ [André (1954)]. As a result, if $p$ is a prime number, there does not exist any nontrivial $(t-1)$-spread. Nevertheless, since we are interested in only $m$ disjoint subspaces, we need only for the maximum number of disjoint $(t-1)$-dimensional subspaces contained in $\mathcal{P}$ to be at least $m$. This is called a partial $(t-1)$-spread in finite projective geometry settings.

DEFINITION 2. A partial $(t-1)$-spread $\mathcal{S}$ of the effect space $\mathcal{P}$ is a set of $(t-1)$-dimensional projective subspaces of $\mathcal{P}$ that are pair-wise disjoint.

Theorem 3 gives necessary and sufficient conditions for the existence of a set $\mathcal{S}$ of disjoint subspaces, each of size $2^t - 1$, for $t < p$. It also provides the size of the minimum overlap in the case where there are no two subspaces that are disjoint. The proof of the theorem can be deduced from that of the more general setup of Theorem 6.

THEOREM 3. *Let $\mathcal{P}$ be a projective space $PG(p-1, 2)$ and let $S_1$ and $S_2$ be two distinct $(t-1)$-dimensional projective subspaces of $\mathcal{P}$ of size $|S_1| = |S_2| = 2^t - 1$ for $0 < t < p$. Then:*

*(a) If $t \leq p/2$, then subspaces $S_1$ and $S_2$ exist such that $S_1 \cap S_2 = \phi$;*

*(b) If $t > p/2$, then for any choice of $S_1, S_2$ in $\mathcal{P}$, $|S_1 \cap S_2| \geq 2^{2t-p} - 1$, and there exist subspaces $S_1, S_2$ such that the equality holds.*

Theorem 3(a) guarantees that, when $t \leq p/2$, one can obtain at least two disjoint $(t-1)$-dimensional subspaces of $\mathcal{P}$. When $p$ is not divisible by $t$, it can be expressed as $p = kt + r$ for nonnegative integers $k, t, r$ satisfying $0 < r < t < p$ and $k \geq 1$. Since $t$ divides $kt$, there must exist a $(t-1)$-spread, $\mathcal{S}_0$, of $PG(kt-1, 2)$ contained in $\mathcal{P}$, where $|\mathcal{S}_0| = (2^{kt} - 1)/(2^t - 1)$. Then, since $\mathcal{S}_0$ is a $(t-1)$-spread of a subspace in $\mathcal{P}$, the set of disjoint $(t-1)$-subspaces in $\mathcal{P}$ may be expandable. The following result ensures the existence of a larger set of disjoint $(t-1)$-dimensional subspaces.

THEOREM 4 [Eisfeld and Storme (2000)]. *Let $\mathcal{P}$ be a finite projective space $PG(p-1, 2)$, with $p = kt + r$ for $0 < r < t < p$. Then, there exists a partial $(t-1)$-spread $\mathcal{S}$ of $\mathcal{P}$ with $|\mathcal{S}| = 2^r \frac{2^{kt}-1}{2^t-1} - 2^r + 1$.*

Next, Theorem 5 summarizes results available on the *maximum number* of pair-wise disjoint $(t-1)$-dimensional subspaces of $\mathcal{P}$ for different combinations of $t$ and $p$.

THEOREM 5 [Govaerts (2005)]. *Let $\mathcal{P}$ be a projective space $PG(p-1, 2)$, with $p = kt + r$ for $k \geq 1$, $0 < r < t < p$, and let $\mathcal{S}$ be a partial $(t-1)$-spread of $\mathcal{P}$ with $|\mathcal{S}| = 2^r \frac{2^{kt}-1}{2^t-1} - s$, where $s$ is known as the deficiency. Then:*
  (a) *If $r = 1$, then $s \geq 2^r - 1 = 1$;*
  (b) *If $r > 1$ and $t \geq 2r$, then $s \geq 2^{r-1} - 1$;*
  (c) *If $r > 1$ and $t < 2r$, then $s \geq 2^{r-1} - 2^{2r-t-1} + 1$.*

Although Theorem 5 provides an upper bound for the number of disjoint subspaces in $\mathcal{P}$, it does not guarantee their existence. As illustrated in the following two examples, the bounds may not be tight.

EXAMPLE 3. Consider a single-replicate $2^5$ factorial experiment with randomization restrictions defined by $S_1, S_2$ and $S_3$ such that $S_1 \supset \{A, B\}$, $S_2 \supset \{C\}$ and $S_3 \supset \{D, E\}$. As discussed in Section 3.2, for analyzing this design, we need at least three half-normal plots depending on the overlapping pattern among the $S_i$'s. The half-normal plot approach requires at least seven effects for each plot [Schoen (1999)]. Therefore, since the $S_i$'s are projective subspaces, the desired features of the RDCSSs are that $S_i \cap S_j = \phi$ for all $i \neq j$ and $|S_i| \geq 7 = 2^3 - 1$ for all $i$. Since $p = 5$ is prime, there does not exist a nontrivial $(t-1)$-spread of $\mathcal{P} = PG(4, 2)$ with $t \geq 3$. However, $p$ can be expressed as $p = kt + r$ with $t = 3$, $k = 1$, $r = 2$ and, from Theorem 5(c), $|\mathcal{S}| \leq 2$. There is no certainty from Theorem 5 regarding the existence of even one pair of disjoint two-dimensional subspaces. In fact, from Theorem 3(b), any pair of two-dimensional subspaces in $PG(4, 2)$ has an overlap of at least $2^{2t-p} - 1 = 1$ effect; hence, the bound is not tight.

Example 4. Consider a $2^8$ full factorial design with $m$ stages of randomization characterized by RDCSSs defined by $S_1, \ldots, S_m$. From Theorem 4, with $p = 8$, $t = 3$, $r = 2$ and $k = 2$, we know that there exists a partial two-spread where the number of disjoint two-dimensional subspaces of the effect space $\mathcal{P}$ is 33. Theorem 5(c) gives the size of maximal partial two-spread of $\mathcal{P}$ to be bounded above by 34. Thus, either the bound in Theorem 5(c) is again not tight or there exists a larger number of disjoint two-spaces of $\mathcal{P}$ than guaranteed by Theorem 4.

Although the results presented in Theorems 3, 4 and 5 focus on $PG(p-1, 2)$, they can easily be generalized to $PG(p-1, q)$, where $q$ is a prime or prime power. See Eisfeld and Storme (2000) and Govaerts (2005) for the generalizations of Theorems 4 and 5, respectively. Theorem 3(a) also holds for two $(t-1)$-dimensional subspaces of $PG(p-1, q)$. For the generalization of Theorem 3(b), if $t > p/2$ and $S_1$, $S_2$ are two $(t-1)$-dimensional subspaces of $PG(p-1, q)$, then $|S_1 \cap S_2| \geq (q^{2t-p} - 1)/(q - 1)$ and there exist a pair such that the equality holds.

We note that, for the special case of 2-level factors, $t = 2$ and $p$ odd, we have $p = 2k + 1$ for positive integer $k$; so, $r = 1$ and Theorem 4 guarantees that the bound of Theorem 5(a) is achieved and $|\mathcal{S}| = (2^p - 5)/3$. For this case, a construction was proposed by Wu (1989) based on the existence of two permutations of the effect space satisfying certain properties. The result of Theorem 4 is more general since it holds for any integer $t$ $(0 < t < p)$ and is easily extendable for arbitrary prime or prime power $q$ in $PG(p-1, q)$.

4.2. *RDCSSs with subspaces of different size.* So far, the results have focused on the existence of disjoint subspaces of the same size. However, it is not unusual for disjoint subspaces of different sizes to be required [see, e.g., Example 4 in Bingham et al. (2008) and the battery cells experiment in Vivacqua, Bisgaard and Steudel (2003)], and we next prove a new result that gives conditions for the existence of such a set of disjoint subspaces.

Theorem 6. *Let $\mathcal{P}$ be a projective space $PG(p-1, 2)$ and $S_i$ be a $(t_i - 1)$-dimensional subspace of $\mathcal{P}$, where $0 < t_i < p$ for $i = 1, 2$. Then:*

(a) *If $t_1 + t_2 \leq p$, then there exist $S_1$ and $S_2$ such that $S_1 \cap S_2 = \phi$;*

(b) *If $t_1 + t_2 > p$, then for any choice of $S_1$, $S_2$ in $\mathcal{P}$, $|S_1 \cap S_2| \geq 2^{t_1+t_2-p} - 1$, and there exist $S_1$ and $S_2$ such that the equality holds.*

Proof. (a) Define the effect space $\mathcal{P} = \langle F_1, \ldots, F_p \rangle$, where $F_i$'s denote the main effects of the independent factors of a $2^p$ factorial experiment. Since $t_1 + t_2 \leq p$, define $S_1 = \langle F_1, \ldots, F_{t_1} \rangle$ and $S_2 = \langle F_{t_1+1}, \ldots, F_{t_1+t_2} \rangle$. Clearly, $S_1$ and $S_2$ are disjoint.

(b) For the case $t_1 + t_2 > p$, define the subspaces $S_1$ and $S_2$ to be $S_1 = \langle F_1, \ldots, F_{t_1} \rangle$ and $S_2 = \langle F_{p-t_2+1}, \ldots, F_{t_1}, F_{t_1+1}, \ldots, F_p \rangle$. Thus, $S_1 \cap S_2 = \langle F_{p-t_2+1}, \ldots, F_{t_1} \rangle$ with $|S_1 \cap S_2| = |PG(t_1 + t_2 - p - 1, 2)| = 2^{t_1+t_2-p} - 1$. Now, if there exist subspaces $S_1^*$ and $S_2^*$ for which $t_1 + t_2 > p$ and $|S_1^* \cap S_2^*| < 2^{t_1+t_2-p} - 1$, then it can be shown that $|\langle S_1, S_2 \rangle| > 2^p - 1$ [see Ranjan (2007), Theorem 4.2 for details]. This contradicts the fact that if $S_1^* \subset \mathcal{P}$ and $S_2^* \subset \mathcal{P}$, then $\langle S_1^*, S_2^* \rangle$ is also a subspace in $\mathcal{P}$. Thus, $S_1$ and $S_2$ provide the minimum possible overlap, as required. □

For $t_1 = t_2 = t$, Theorem 6 simplifies to Theorem 3. When $t_1 + t_2 \leq p$ [as in Theorem 6(a)], in addition to $S_1$ and $S_2$, one can expect mutually disjoint subspaces of size $2^t - 1$ that do not overlap with $S_1$ and $S_2$, where $t < p - \max(t_1, t_2)$. The next theorem, which is the main result of this section, establishes the existence of a set of unequal sized subspaces of $\mathcal{P}$, where $t_j + t_k \leq p$ for any pair of subspaces $S_j$ and $S_k$.

THEOREM 7. *Let $\mathcal{P}$ be a projective space $PG(p-1, 2)$ and $S_1$ be a $(t_1-1)$-dimensional subspace of $\mathcal{P}$. If $p/2 < t_1 < p$, then there exist $S_2, \ldots, S_m$ such that $S_i \cap S_j = \phi$ for all $i, j \in \{1, \ldots, m\}$, where $|S_i| = 2^{t_i} - 1$ for $t_i \leq p - t_1$, $2 \leq i \leq m$ and $m = 2^{t_1} + 1$.*

PROOF. Let $p/2 < t_1 < p$ and define $t^* = p - t_1$. Then, the effect space $\mathcal{P}$ is a $PG(t_1 + t^* - 1, 2)$. Let $\mathcal{P}' = PG(2t_1 - 1, 2)$ be such that $\mathcal{P}' \supseteq \mathcal{P}$. Let $S_1$ be a $(t_1-1)$-dimensional subspace of $\mathcal{P}$, and let $\mathcal{S}'$ be a $(t_1-1)$-spread of $\mathcal{P}'$ that contains $S_1$. Then,

$$|\mathcal{S}'| = m = \frac{|PG(2t_1 - 1, 2)|}{|PG(t_1 - 1, 2)|} = \frac{2^{2t_1} - 1}{2^{t_1} - 1} = 2^{t_1} + 1.$$

The set of disjoint $(t^* - 1)$-spaces of $\mathcal{P}$ whose elements are disjoint from $S_1$ is given by $\mathcal{S}^* = \{S \cap \mathcal{P} : S \in \mathcal{S}' \setminus \{S_1\}\}$. The elements of $\mathcal{S}^*$ can be denoted by $S_2^*, \ldots, S_m^*$. For every $i = 2, \ldots, m$, since the desired subspace $S_i$ is of size $2^{t_i} - 1$ with $t_i \leq p - t_1 = t^*$, one can construct $S_i$ by constructing a $(t_i - 1)$-dimensional subspace of $S_i^*$. □

This theorem is important as it guarantees the existence of $2^{t_1} + 1$ disjoint subspaces of different sizes, with one $(t_1-1)$-dimensional subspace ($p/2 < t_1 < p$) and $2^{t_1}$ subspaces of dimension up to $(t-1)$ each, where $t \leq p - t_1$. Note that the boundary conditions $t_1 = p$ and $t_1 = p/2$ are uninteresting, as $t_1 = p$ implies that $S_1$ is the entire effect space, and $t_1 = p/2$ guarantees the existence of a $(t_1-1)$-spread of $\mathcal{P}$. The proof of Theorem 7 leads to a construction strategy for $m = 2^{t_1} + 1$ disjoint subspaces of unequal sizes, as illustrated in Section 5.3. Although Theorem 4 is not a special case of Theorem 7, a similar construction is applicable for equal sized RDCSSs and is discussed in Section 5.2.

Both Theorems 6 and 7 can be generalized to $q$-level designs for prime or prime power $q$. Theorem 6(a) also holds for $(t_i - 1)$-dimensional subspaces $S_i$'s in $PG(p-1, q)$. For Theorem 6(b), if $t_1 + t_2 > p$, then $|S_1 \cap S_2| \geq (q^{t_1+t_2-p}-1)/(q-1)$, and there exists a pair $S_1, S_2$ such that the equality holds. The generalization of Theorem 7 guarantees the existence of a set of $m = q^{t_1} + 1$ disjoint subspaces of different sizes, where $|S_i| = (q^{t_i}-1)/(q-1)$ for $p/2 < t_1 < p$ and $t_i \leq p - t_1$. The proofs of results for arbitrary prime or prime power $q$ are similar to the $q = 2$ case shown above.

Thus far, we have established some necessary and some sufficient conditions for the existence of a set of disjoint subspaces of the same size and also of different sizes. If the desired number of stages of randomization is less than or equal to the number of subspaces guaranteed to exist, we can obtain an appropriate subset satisfying the randomization restrictions required by the experimenter. Finding an actual design with the required properties is the next issue for the experimenter, and this is discussed in the next section.

**5. Construction of disjoint subspaces.** A construction approach is now proposed for factorial designs with $m$ levels of randomization. First, the construction for equal sized disjoint subspaces is presented, followed by the construction of disjoint subspaces of different sizes. The subspaces themselves have no statistical meaning until the factors have been assigned to columns of the model matrix or, equivalently, to points in $PG(p-1, 2)$. The set of disjoint subspaces obtained from an arbitrary assignment may not directly satisfy the experimenter's requirements for the factor levels. Consequently, we propose an algorithm that transforms a set of subspaces obtained from the RDCSS construction to a set of subspaces satisfying the desired restrictions on the factor levels.

5.1. *RDCSSs and $(t-1)$-spreads.* When $t$ divides $p$, we know that there exists a $(t-1)$-spread of $\mathcal{P} = PG(p-1, 2)$ [André (1954)]. The construction of the spread starts with writing the $2^p - 1$ nonzero elements of $GF(2^p)$ in cycles of length $N$ [Hirschfeld (1998)]. An element $w$ is called primitive if $\{w^i : i = 0, 1, \ldots, u - 2\} = GF(u) \setminus \{0\}$. Let $w$ be a root of a primitive polynomial of degree $p$ for $GF(2^p)$. Then, the $2^p - 1$ elements of the effect space $\mathcal{P}$ or, equivalently, the nonzero elements of $GF(2^p)$ are $w^i$, $i = 0, \ldots, 2^p - 2$. The element $w^i$ can be written as a linear combination of the basis monomials $w^0, \ldots, w^{p-1}$. The element

$$w^i = \alpha_0 w^{p-1} + \alpha_1 w^{p-2} + \cdots + \alpha_{p-2} w + \alpha_{p-1} \tag{7}$$

represents an $r$-factor interaction $\delta = (\alpha_0, \alpha_1, \ldots, \alpha_{p-1})$ for $\alpha_i \in GF(2)$ if exactly $r$ elements in $\delta$ are nonzero. Following this representation for the factorial effects in $\mathcal{P}$, and using shorthand notation $(iN + j)$ to denote $w^{iN+j}$ for $0 \le i \le \theta - 1$, $0 \le j \le N - 1$ (where $\theta$ is the number of cycles), cycles of length $N$ can be written as in Table 1. The columns define the projective subspaces $S_1, \ldots, S_N$. A necessary and sufficient condition [Hirschfeld (1998), Chapter 4] that there exists a $(t-1)$-space of cycles of length $N$ which is smaller than $|PG(p-1,2)|$ is that the greatest common divisor, $gcd(t,p)$, of $t$ and $p$ is greater than one. Then, $N = |PG(p-1,2)|/|PG(l-1,2)|$, where $l = gcd(t,p)$. Thus, when $t$ divides $p$, there exist $2^t - 1$ cycles each of length $N$ which lead to the formation of the required subspaces $S_1, \ldots, S_N$ of $\mathcal{P}$ with $S_i \cap S_j = \phi$ for $i \ne j$; that is, $S_1, \ldots, S_N$ constitute a $(t-1)$-spread $\mathcal{S}$ of the effect space $\mathcal{P} = PG(p-1,2)$.

TABLE 1
*Cycles of length $N$*

| $S_1$ | $S_2$ | ...... | $S_N$ |
|---|---|---|---|
| 0 | 1 | ...... | $N-1$ |
| $N$ | $N+1$ | ...... | $2N-1$ |
| $\vdots$ | $\vdots$ | $\vdots$ | $\vdots$ |
| $(\theta-1)N$ | $(\theta-1)N+1$ | ...... | $\theta N - 1$ |

A $(t-1)$-spread of $PG(p-1,2)$, obtained as above, distributes all the factor main effects evenly among the $|\mathcal{S}|$ disjoint subspaces. However, restrictions on the $m$ stages of randomization are usually pre-specified by the experimenter and, as illustrated in Example 5, an RDCSS for a block design will contain no main effects whereas, for a split-lot design, several factor main effects may be assigned to one or more RDCSS.

EXAMPLE 5. Consider a single-replicate $2^6$ factorial experiment with the randomization structure determined by a *blocked split-lot design*, where the experiment has to be performed in blocks of size eight each. Suppose that the experimenter wishes to specify the factorial effects, $ABC, BDE$ and $CEF$, to be confounded with the blocks (i.e., $S_1^* = \langle ABC, BDE, CEF \rangle$). In addition, suppose the trials proceed in a two-step process, where the restrictions imposed by the experimenter on the two steps of randomization are such that $S_2^* \supset \{A, B\}$ and $S_3^* \supset \{D\}$. As a result, there are three restrictions on the randomization of the experiment: one due to blocking the experimental units ($S_1^*$) and the other two due to splitting the experimental units into sub-lots ($S_2^*, S_3^*$). To use half-normal plots, it is desirable to have three disjoint subspaces each of at least size seven (i.e., $t = 3$), where the subspaces should satisfy the restrictions defined by $S_1^*, S_2^*$ and $S_3^*$. Here $gcd(t, p) = 3$, so $N = |PG(5, 2)|/|PG(2, 2)| = (2^6 - 1)/(2^3 - 1) = 9$, and there exist 7 cycles of length 9 or, equivalently, 9 disjoint subspaces of size 7 (i.e., a two-spread of $\mathcal{P}$). The two-spread $\mathcal{S} = \{S_1, \ldots, S_9\}$ obtained using the primitive polynomial $w^6 + w + 1$, root $w$ and representation (7) is shown in Table 2. Notice that each two-dimensional subspace of $\mathcal{P}$ in $\mathcal{S}$ contains at most one main effect, and none matches the design requirements. It is shown below how to transform the spread $\mathcal{S}$ to $\mathcal{S}^*$, such that $\mathcal{S}^*$ contains the three disjoint subspaces $S_1^*, S_2^*$ and $S_3^*$ satisfying the experimenter's requirement.

TABLE 2
*The two-spread, $\mathcal{S}$, of $PG(5,2)$ using the primitive polynomial $w^6 + w + 1$*

| $S_1$ | $S_2$ | $S_3$ | $S_4$ | $S_5$ | $S_6$ | $S_7$ | $S_8$ | $S_9$ |
|-------|-------|-------|-------|-------|-------|-------|-------|-------|
| F     | E     | D     | C     | B     | A     | EF    | DE    | CD    |
| BC    | AB    | AEF   | DF    | CE    | BD    | AC    | BEF   | ADE   |
| CDEF  | BCDE  | ABCD  | ABCEF | ABDF  | ACF   | BF    | AE    | DEF   |
| CDE   | BCD   | ABC   | ABEF  | ADF   | CF    | BE    | AD    | CEF   |
| BDE   | ACD   | BCEF  | ABDE  | ACDEF | BCDF  | ABCE  | ABDEF | ACDF  |
| BCF   | ABE   | ADEF  | CDF   | BCE   | ABD   | ACEF  | BDF   | ACE   |
| BDEF  | ACDE  | BCDEF | ABCDE | ABCDEF| ABCDF | ABCF  | ABF   | AF    |

For the transformation of spreads, we use an appropriate *collineation* [see, e.g., Batten (1997)] of the projective space $\mathcal{P}$. A collineation of $PG(p-1, q)$ is a permutation $f$ of its points such that $(t-1)$-dimensional subspaces are mapped to $(t-1)$-dimensional subspaces for $1 \le t \le p$. The existence of a collineation $f$ from $\mathcal{S}$ to $\mathcal{S}^*$ is equivalent to the existence of a $p \times p$ matrix $\mathcal{M}$ such that, for every given $S \in \mathcal{S}$, there is a unique $S^* \in \mathcal{S}^*$ and, for every $z \in S$, there exists a unique $z^* \in S^*$ that satisfies $z^{*\prime} = z'\mathcal{M}$. Note that the transformation of a spread amounts to relabelling the columns of the model matrix. As a result, one cannot find a collineation matrix $\mathcal{M}$ if the experimenter's requirement is not feasible. Moreover, if the desired set of subspaces is nonisomorphic to the spread we started with, then also there does not exist any relabelling to obtain the desired design. However, finding an appropriate collineation matrix, whenever it exists, is also nontrivial. Next, we propose an algorithm that finds a collineation matrix $\mathcal{M}$, if it exists, and concludes the nonexistence if one does not exist.

The proposed algorithm is illustrated through the setup of the $2^6$ experiment of Example 5. To obtain a set of disjoint subspaces satisfying the restrictions imposed on the three stages of randomization, we have to find an appropriate $6 \times 6$ collineation matrix $\mathcal{M}$. The proposed algorithm for finding the matrix $\mathcal{M}$ is outlined as follows:

1. Select one of the $\binom{9}{3}$ possible choices from the spread $\mathcal{S}$ in Table 2 to be a set of three disjoint subspaces. For example, suppose that $S_1, S_3$ and $S_7$ are chosen for $S_2^*, S_3^*$ and $S_1^*$, respectively.

2. Choose two effects from $S_1$, one effect from $S_3$ and three effects from $S_7$ to relabel these to the desired effects $(A, B)$, $D$ and $(ABC, BDE, CEF)$ in $S_2^*, S_3^*$ and $S_1^*$, respectively. For example, one choice among the $\binom{7}{2}\binom{7}{1}\binom{7}{3}$ different options is $\{CDE, BCF, D, EF, AC, BF\}$. The collineation matrix is defined by the mapping induced by $CDE \to A, BCF \to B$, $D \to D, \ldots, BF \to CEF$.

3. Denote the $(i, j)$th entry of the $p \times p$ matrix $\mathcal{M}$ as $x_k$ for $k = j + (i - 1)p$ (i.e., list the elements row by row). Construct a $p^2 \times p^2$ matrix $\mathcal{Q}$ and a $p^2 \times 1$ vector $\delta$ as follows. Define the rows of matrix $\mathcal{Q}$ and vector $\delta$ in the order of restrictions on the transformation. For the example under consideration, the first transformation, $(CDE)'\mathcal{M} = (A)'$, can be written as

$$(8) \qquad \begin{bmatrix} 0 & 0 & 1 & 1 & 1 & 0 \end{bmatrix} \mathcal{M} = \begin{bmatrix} 1 & 0 & 0 & 0 & 0 & 0 \end{bmatrix}.$$

In total, there are $p$ $(1 \leq s \leq p)$ independent transformations, and for each transformation there are $p$ $(p(s - 1) \leq i \leq ps)$ rows of $\delta$ and $\mathcal{Q}$. The first set $(s = 1)$ of $p = 6$ elements of $\delta$ are given by the right-hand side of (8). The corresponding rows of $\mathcal{Q}$ are defined as

$$\mathcal{Q}_{il} = 1, \qquad \text{if } l = (\tau - 1)p + \text{mod}(i - 1, p) + 1$$
$$\text{and the } \tau\text{th entry of } \begin{bmatrix} 0 & 0 & 1 & 1 & 1 & 0 \end{bmatrix} \text{ is nonzero,}$$
$$= 0, \qquad \text{otherwise,}$$

for $1 \leq i \leq p$. Similarly, all the rows of the matrix $\mathcal{Q}$ and the vector $\delta$ can be expressed using $p$ restrictions on the transformation. For example, the second set $(s = 2)$ of $p = 6$ rows $(\delta_i, \mathcal{Q}_{il}$ for $p + 1 \leq i \leq 2p)$ are given by the right-hand side of $(BCF)'\mathcal{M} = (B)'$; that is, $\begin{bmatrix} 0 & 1 & 1 & 0 & 0 & 1 \end{bmatrix} \mathcal{M} = \begin{bmatrix} 0 & 1 & 0 & 0 & 0 & 0 \end{bmatrix}$, and so on for the remaining transformations.

4. If there exists a solution for $\mathcal{Q}x = \delta$, reconstruct the $p \times p$ matrix $\mathcal{M}$ (row by row) from the $p^2 \times 1$ solution vector $x = \mathcal{Q}^{-L}\delta$, where $\mathcal{Q}^{-L}$ is a left inverse of $\mathcal{Q}$. Exit the algorithm.

5. If there does not exist a solution then go to step 2 and, if possible, choose a different set of effects from the subspaces selected in step 1.

6. If all possible choices for the set of effects from these three subspaces have been exhausted, go to step 1 and choose a different set of three subspaces.

7. If all the $\binom{9}{3}$ choices for a set of effects have been used and still a solution does not exist, then either the two spreads $\mathcal{S}$ and $\mathcal{S}^*$ are nonisomorphic or the experimenter's requirement is not achievable. Thus, the desired spread cannot be obtained from $\mathcal{S}$.

In the illustration used above, the effects chosen for relabelling the columns to achieve the desired design do provide a feasible solution to $\mathcal{Q}x = \delta$. The collineation matrix $\mathcal{M}$, obtained from the solution $x = \mathcal{Q}^{-L}\delta$, is given by

$$\mathcal{M} = \begin{pmatrix} 0 & 0 & 1 & 1 & 0 & 1 \\ 0 & 0 & 1 & 1 & 0 & 0 \\ 0 & 1 & 1 & 0 & 1 & 1 \\ 0 & 0 & 0 & 1 & 0 & 0 \\ 1 & 1 & 1 & 1 & 1 & 1 \\ 0 & 0 & 0 & 1 & 1 & 1 \end{pmatrix}.$$

In this example, an exhaustive search found that 45.7% of all possible choices give a feasible solution to the equation $\mathcal{Q}x = \delta$. That is, an arbitrary choice of $p$ independent effects from $\mathcal{S}$ (according to steps 1 and 2) results in a feasible design 45.7% of the time. The rest of the time, the solution is infeasible since the full factorial design becomes a replicated fraction.

The spread acts as a template which enables a faster search than the exhaustive relabelling of all the factorial effects to find the design satisfying the experimenter's requirement. For this example, our algorithm requires *at most* $\binom{9}{3}\binom{7}{2}\binom{7}{1}\binom{7}{3} \approx 5 \times 10^5$ different relabellings, whereas an exhaustive relabelling approach can require up to $(2^6 - 1)! \approx 2 \times 10^{87}$ different relabellings. While searching through each of the relabellings in the template can be time consuming, our Matlab 7.0.4 implementation of the algorithm found the first feasible collineation matrix in 5.34 seconds on a Pentium(R) 4 processor machine running Windows XP.

5.2. *Partial $(t-1)$-spreads.* When $t$ does not divide $p$, Theorem 4 guarantees the existence of $|\mathcal{S}| = 2^r \frac{2^{kt}-1}{2^t-1} - 2^r + 1$ disjoint subspaces, where $p = kt + r$. For constructing these subspaces, one can use the existence proof of Eisfeld and Storme (2000) for the most part. However, the proof assumes the existence of an $(s_i - 1)$-spread $\mathcal{S}'_i$ of $\mathcal{P}_i$ that contains an $(s_i - 1)$-dimensional subspace, $U_i$, of $\mathcal{P}'_{i+1}$, where $s_i = it + r$, $\mathcal{P}_i = PG(2s_i - 1, 2)$ and $\mathcal{P}'_{i+1} = PG(s_i + t - 1, 2)$ for $i = 1, \ldots, k-1$. This is nontrivial, and we propose a two-step construction method to get this as follows: (a) construct an $(s_i - 1)$-spread $\mathcal{S}''_i$ of $\mathcal{P}_i$ as described in Section 5.1 and then (b) transform the spread $\mathcal{S}''_i$ to $\mathcal{S}'_i$ by finding an appropriate collineation such that $U_i \in \mathcal{S}'_i$. Thus, we can construct a set of $|\mathcal{S}|$ disjoint $(t-1)$-dimensional subspaces. Finally, we find an appropriate collineation to transform the partial $(t-1)$-spread $\mathcal{S}$ to obtain the $m$ RDCSSs satisfying the experimenter's requirement.

5.3. *Disjoint subspaces of different sizes.* A more general setting is when the RDCSSs are allowed to have different sizes. For a $2^p$ full factorial design, Theorem 7 guarantees the existence of one subspace $S_1$ of size $2^{t_1} - 1$ with $t_1 > p/2$ and $2^{t_1}$ subspaces of size at most $2^{p-t_1} - 1$. For constructing these $2^{t_1} + 1$ mutually disjoint subspaces of $\mathcal{P}$, the proof of Theorem 7 requires constructing a $(t_1 - 1)$-spread $\mathcal{S}''$ of $PG(2t_1 - 1, 2)$ that contains $S_1$. The spread $\mathcal{S}''$ can be obtained by first constructing a $(t_1 - 1)$-spread of $PG(2t_1 - 1, 2)$ and then by applying the appropriate collineation ($\mathcal{M}_0$) found by the algorithm described in Section 5.1 to obtain $\mathcal{S}'$ that contains $S_1$. After $\mathcal{S} = \{S \cap \mathcal{P} : S \in \mathcal{S}' \setminus \{S_1\}\}$ is obtained, one has to find a suitable collineation ($\mathcal{M}_1$) so that the final set of subspaces satisfy the experimenter's restrictions on RDCSSs. The algorithm can be made more efficient by combining the problem of finding the two collineation matrices into one. When transforming the spread $\mathcal{S}''$ to $\mathcal{S}'$ containing $S_1$, we can impose other restrictions on $S_2, \ldots, S_m$ in this step itself. The steps of the construction are illustrated through Example 6.

TABLE 3
*Three-spread of $PG(7,2)$ containing $S_1 = \langle A,B,C,D \rangle$, $S_2 \supset \{E,F\}$ and $S_3 \supset \{G\}$*

| $S_1$ | $S_2$ | $S_3$ | $S_4$ | $\cdots$ | $S_{16}$ | $S_{17}$ |
|-------|-------|-------|-------|----------|----------|----------|
| A     | E     | BCFG  | DE    | $\cdots$ | DEFH     | CE       |
| B     | AEFH  | G     | BEFH  | $\cdots$ | ACFG     | ADEFH    |
| C     | CFG   | ADFH  | BFG   | $\cdots$ | BCG      | ABCFG    |
| D     | CG    | BCEGH | DG    | $\cdots$ | ACFH     | BG       |
| AB    | AFH   | BCF   | BDFH  | $\cdots$ | ACDEGH   | ACDFH    |
| BC    | ACEGH | ADFGH | EGH   | $\cdots$ | ABF      | BCDEGH   |
| CD    | F     | ABCDEFG | BDF | $\cdots$ | ABFGH    | ACF      |
| ABD   | ACFGH | EFGH  | BFGH  | $\cdots$ | DEFG     | ABCDFGH  |
| AC    | CEFG  | ABCDGH | BDEFG | $\cdots$ | BCDEFGH  | ABEFG    |
| BD    | ACEFGH | BCEH | BDEFGH | $\cdots$ | GH       | ABDEFGH  |
| ABC   | ACGH  | ABCDH | DGH   | $\cdots$ | ABDEH    | BDGH     |
| BCD   | AEH   | ABCDEF | DEH  | $\cdots$ | BCH      | CDEH     |
| ABCD  | AH    | ADEG  | H     | $\cdots$ | BCDEF    | DH       |
| ACD   | EF    | ADE   | BEF   | $\cdots$ | ABDEG    | AEF      |
| AD    | CEG   | EFH   | EG    | $\cdots$ | ACDE     | BCEG     |

EXAMPLE 6.    Consider a single-replicate $2^7$ factorial design with three stages of randomization. Let the restrictions imposed by the experimenter on the three RDCSSs be $S_1 \supset \{A,B,C,D\}$, $S_2 \supset \{E,F\}$ and $S_3 \supset \{G\}$. Following the notation of Theorem 7, since $p = 7$ and $t_1 = 4$, there exist $m = 2^4 + 1 = 17$ pair-wise disjoint subspaces with $|S_i| = 2^{t_i} - 1$ for $i = 1, \ldots, 17$, where $t_1 = 4$ and $t_i \leq 3$ for $i = 2, \ldots, 17$. First, we construct a three-spread $\mathcal{S}''$ of $PG(7,2)$ using the method described in Section 5.1, and then we find an appropriate collineation matrix $\mathcal{M}_0$ that transforms $\mathcal{S}''$ to $\mathcal{S}$ such that $\mathcal{S}$ contains $S_1 = \langle A,B,C,D \rangle$, $S_2 \supset \{E,F\}$ and $S_3 \supset \{G\}$. Table 3 contains some of the elements of the transformed spread.

For this example, the three disjoint subspaces $S_1, S_2$ and $S_3$ that satisfy the experimenter's requirements are obtained by deleting the elements of $S_i$'s in the transformed spread that contains $H$. The collineation matrix $\mathcal{M}$ used for the transformation is as follows:

$$\mathcal{M} = \begin{pmatrix} 1 & 1 & 1 & 1 & 1 & 1 & 1 & 0 \\ 0 & 1 & 0 & 0 & 1 & 1 & 0 & 0 \\ 1 & 0 & 1 & 1 & 1 & 0 & 1 & 0 \\ 1 & 1 & 1 & 0 & 1 & 0 & 0 & 1 \\ 0 & 1 & 0 & 1 & 0 & 1 & 1 & 1 \\ 0 & 1 & 1 & 1 & 1 & 1 & 0 & 0 \\ 0 & 1 & 0 & 0 & 0 & 1 & 1 & 1 \\ 1 & 0 & 0 & 0 & 0 & 0 & 0 & 0 \end{pmatrix}.$$

Thus, $\{S \cap \mathcal{P} : S \in \mathcal{S}' \setminus \{S_1\}\} \cup S_1$ contains the required set of subspaces $S_1, S_2$ and $S_3$ for the three stages of randomization.

**6. Fractional factorial designs.** Though the designs discussed in this article are full factorial designs, the results also apply to FF designs with randomization restrictions. The existence of regular FF designs with randomization restrictions is equivalent to that of the full factorial design generated from the basic factors of the FF design (*base factorial design*). One can find a $2^{r-s}$ fractional factorial with randomization restrictions by first investigating the existence a $2^u$ full factorial design ($u = r - s$) with the appropriate randomization restrictions. If the design exists, one can construct the $2^u$ full factorial design with randomization restrictions as outlined in Section 5 and then assign the *added factors* to the effects at the appropriate level of randomization to get the desired design.

Consider, for example, a $2^{8-2}$ fractional factorial experiment with randomization structure characterized by a split-lot design. Suppose that the experiment has to be run in four stages with randomization restrictions given by $S_1 \supset \{A, B\}$, $S_2 \supset \{C, D\}$, $S_3 \supset \{E, F\}$ and $S_4 \supset \{G, H\}$. Then, the six independent *basic factors* $(A, B, C, D, E, F)$ with effect space $\mathcal{P} = \langle A, B, \ldots, F \rangle$ form a $2^6$ full factorial split-lot design. The results discussed in Section 4 guarantee the existence of a two-spread of $\mathcal{P}$, and the algorithm described in Section 5.1 shows the construction for obtaining three disjoint subspaces of size seven each satisfying the requirements of $S_1$, $S_2$ and $S_3$. Moreover, we know that there exist nine disjoint subspaces of size seven each. Therefore, $S_4$ can be constructed by choosing a subspace from the remaining six disjoint subspaces and then by aliasing two generators $G$ and $H$ with effects in this subspace [Box, Hunter and Hunter (1978)].

Fractionation of the base factorial design can occur in a different scenario. For example, consider a $2^{8-2}$ experiment with the requirement of three stages of randomization and the restrictions on the RDCSSs be defined by $S_1 \supset \{A, B\}$, $S_2 \supset \{C, D, E\}$ and $S_3 \supset \{F, G, H\}$. In this case, one can use the algorithm in Section 5.1 to construct three disjoint subspaces satisfying $S_1 \supset$ $\{A, B\}$, $S_2 \supset \{C, D, E\}$ and $S_3 \supset \{F\}$. Next, one can choose two generators (or points) from $S_3$ that are assigned to $G$ and $H$.

To rank the designs, one can use existing criteria such as minimum aberration [Fries and Hunter (1980)], maximum number of clear effects [Chen, Sun and Wu (1993) and Wu and Chen (1992)] and $V$-criterion [Bingham et al. (2008)]. Thus, one can select an appropriate set of generators based on the experimenter's interest. The structure of a spread acts as a template to shorten the computer search for good fractional factorial designs.

**7. Concluding remarks.** In this paper, we have demonstrated that the projective subspaces of the effect space $\mathcal{P}$ can be used to characterize the randomization restrictions of factorial designs in block, split-plot, split-lot and other structures. Under the assumptions of model (1), Theorem 2 summarizes the impact of randomization restrictions on the distribution of factorial effect estimators and motivates the search for disjoint randomization defining contrast subspaces. Obtaining a set of disjoint subspaces of the effect space $\mathcal{P}$ is nontrivial. In the most general case, Theorem 7 presents a necessary and sufficient condition for the existence of disjoint subspaces of unequal sizes. Furthermore, the proof motivates strategies for constructing these, and this makes a wide variety of designs accessible in this unified framework.

**Acknowledgments.** We would like to thank the referees and Associate Editor for their thoughtful comments and suggestions.

## REFERENCES

ADDELMAN, S. (1964). Some two-level factorial plans with split-plot confounding. *Technometrics* **6** 253–258.

ANDRÉ, J. (1954). Uber nicht-Desarguessche Ebenen mit transitiver Translationsgruppe. *Math. Z.* **60** 156–186. MR0063056

BAILEY, R. A. (1977). Patterns of confounding in factorial designs. *Biometrika* **64** 597–604. MR0501643

BATTEN, L. M. (1997). *Combinatorics of Finite Geometries*, 2nd ed. Cambridge Univ. Press, Cambridge. MR1474497

BINGHAM, D., SITTER, R. R., KELLY, E., MOORE, L. and OLIVAS, J. D. (2008). Factorial designs with multiple levels of randomization. *Statist. Sinica* **18** 493–513.

BINGHAM, D. and SITTER, R. R. (1999). Minimum-aberration two-level fractional factorial split-plot designs. *Technometrics* **41** 62–70.

BISGAARD, S. (1994). A note on the definition of resolution for blocked $2^{k-p}$ designs. *Technometrics* **36** 308–311. MR1292444

BISGAARD, S. (2000). The design and analysis of $2^{k-p} \times 2^{q-r}$ split-plot experiments. *J. Quality Technology* **32** 39–56.

BOSE, R. C. (1947). Mathematical theory of the symmetrical factorial design. *Sankhyā* **8** 107–166.

Box, G. E. P., Hunter, W. G. and Hunter, J. S. (1978). *Statistics for Experiments: An Introduction to Design, Data Analysis, and Model Building.* Wiley, New York. MR0483116

Butler, N. A. (2004). Construction of two-level split-lot fractional factorial designs for multistage processes. *Technometrics* **46** 445–451. MR2101512

Chen, H. and Cheng, C. S. (1999). Theory of optimal blocking of $2^{n-m}$ designs. *Ann. Statist.* **27** 1948–1973. MR1765624

Chen, J., Sun, D. X. and Wu, C. F. J. (1993). A catalogue of two-level and three-level fractional factorial designs with small runs. *Internat. Statist. Rev.* **61** 131–145.

Cheng, S.-W., Li, W. and Ye, K. Q. (2004). Blocked nonregular two-level factorial designs. *Technometrics* **46** 269–279. MR2082497

Cochran, W. G. and Cox, G. M. (1957). *Experimental Designs*, 2nd ed. Wiley, New York. MR0085682

Daniel, C. (1959). Use of half normal plots in interpreting factorial two-level experiments. *Technometrics* **1** 311–341. MR0125710

Dean, A. M. (1978). The analysis of interactions in single replicate generalized cyclic designs. *J. Roy. Statist. Soc. Ser. B* **40** 79–84. MR0512146

Dean, A. M. and John, J. A. (1975). Single-replicate factorial experiments in generalized cyclic designs, II. Asymmetrical arrangements. *J. Roy. Statist. Soc. Ser. B* **37** 72–76. MR0368342

Dean, A. M. and Voss, D. T. (1999). *Design and Analysis of Experiments.* Springer, New York. MR1673800

Dey, A. and Mukerjee, R. (1999). *Fractional Factorial Plans.* Wiley, New York. MR1679441

Eisfeld, J. and Storme, L. (2000). (Partial) $t$-spreads and minimal $t$-covers in finite projective spaces. Univ. Gent.

Fries, A. and Hunter, W. G. (1980). Minimum aberration $2^{k-p}$ designs. *Technometrics* **22** 601–608. MR0596803

Govaerts, P. (2005). Small maximal partial $t$-spreads. *Bull. Belg. Math. Soc. Simon Stevin* **12** 607–615. MR2206003

Hirschfeld, J. W. P. (1998). *Projective Geometries Over Finite Fields.* Oxford Univ. Press, Oxford. MR1612570

Huang, P., Chen, D. and Voelkel, J. O. (1998). Minimum-aberration two-level split-plot designs. *Technometrics* **40** 314–326. MR1659353

Mee, R. W. and Bates, R. L. (1998). Split-lot designs: Experiments for multistage batch processes. *Technometrics* **40** 127–140.

Miller, A. (1997). Strip-plot configurations of fractional-factorials. *Technometrics* **39** 153–161. MR1452344

Milliken, G. A. and Johnson, D. E. (1984). *Analysis of Messy Data—Designed Experiments* **1**. Chapman & Hall, New York.

Mukerjee, R. and Wu, C. F. J. (1999). Blocking in regular fractional factorials: A projective geometric approach. *Ann. Statist.* **27** 1256–1271. MR1740111

Ranjan, P. (2007). Factorial and fractional factorial designs with randomization restrictions—a projective geometric approach. Ph.D. thesis, Simon Fraser Univ.

Schoen, E. D. (1999). Designing fractional two-level experiments with nested error structures. *J. Appl. Statist.* **26** 495–508.

Sitter, R. R., Chen, J. and Feder, M. (1997). Fractional resolution and minimum aberration in blocked factorial designs. *Technometrics* **39** 382–390. MR1482516

Srivastava, J. N. (1987). Advances in the general theory of factorial designs based on partial pencils in Euclidean $n$-space. *Utilitas Math.* **32** 75–94. MR0921637

STAPLETON, R. D., LEWIS, S. M. and DEAN, A. M. (2009). Construction of fractional factorial experiments in split-plot row-column designs. Submitted for publication.

SUN, D. X., WU, C. F. J. and CHEN, Y. Y. (1997). Optimal blocking schemes for $2^n$ and $2^{n-p}$ designs. *Technometrics* **39** 298–307. MR1462588

TAGUCHI, G. (1987). *System of Experimental Design: Engineering Methods to Optimize Quality and Minimize Costs* **1–2**. UNIPUB/Kraus International Pub., White Plains, NY.

VIVACQUA, C. A., BISGAARD, S. and STEUDEL, H. J. (2003). Post-fractionated strip-block designs: A tool for robustness applications and multistage processes. In *2003 Quality & Prod. Research Conference—IBM T. J. Watson Research Ctr.* Yorktown Heights, NY, USA.

VOSS, D. T. and DEAN, A. M. (1987). A comparison of classes of single replicate factorial designs. *Ann. Statist.* **15** 376–384. MR0885743

WU, C. F. J. (1989). Construction of $2^m 4^n$ designs via a grouping scheme. *Ann. Statist.* **17** 1880–1885. MR1026317

WU, C. F. J. and CHEN, Y. (1992). A graph-aided method for planning two-level experiments when certain interactions are important. *Technometrics* **34** 162–175.

P. RANJAN
DEPARTMENT OF MATHEMATICS
  AND STATISTICS
ACADIA UNIVERSITY
WOLFVILLE
NOVA SCOTIA B4P 2R6
CANADA
E-MAIL: pritam.ranjan@acadiau.ca

D. R. BINGHAM
DEPARTMENT OF STATISTICS
  AND ACTUARIAL SCIENCE
8888 UNIVERSITY DRIVE
BURNABY BC V5A 1S6
CANADA
E-MAIL: dbingham@stat.sfu.ca

A. M. DEAN
DEPARTMENT OF STATISTICS
OHIO STATE UNIVERSITY
1958 NEIL AVENUE
COLUMBUS, OHIO 43210
USA
E-MAIL: dean.9@osu.edu